\documentclass[11pt]{amsart}
\usepackage{amsmath, amssymb}
\usepackage{amsfonts}
\usepackage{amsthm}
\usepackage{amsopn}
\usepackage{amscd}

\usepackage [all]{xy}

\newcommand{\sC}{\mathcal{C}}
\newcommand{\cinf}{\sC^\infty}

\newcommand{\sP}{\mathcal{P}}

\newcommand{\bR}{\mathbf{R}}

\newcommand{\comp}{\circ}

\newcommand{\ti}{\tilde}

\newcommand{\diff}{\text{Diff}}

\theoremstyle{plain}
\newtheorem{thm}{Theorem}[section]
\newtheorem{prop}[thm]{Proposition}

\newtheorem{cor}[thm]{Corollary}

\theoremstyle{definition}

\theoremstyle{remark}

\let\a\alpha
\let\b\beta
\let\g\gamma

\let\f\varphi

\let\G\Gamma
\let\D\Delta

\let\na\nabla
\let\ra\rightarrow

\let\ti\tilde

\begin{document}
\title{Another look at connections}
\author{Florin Dumitrescu}
\date{\today}
\maketitle

\begin{abstract} In this note we make use of some properties of vector fields on a manifold to give an alternate proof to \cite{D4} for the equivalence between connections and parallel transport on vector bundles over manifolds. Out of the proof will emerge a new approach to connections on a bundle as a consistent way to lift the dynamics of the manifold to the bundle.
\end{abstract}

\vspace{.2in}

This note is aimed at proving the equivalence between connections (covariant derivatives) and the geometric notion of parallel transport. A classical proof of this fact appears in \cite{D4}. A functorial approach to parallel transport and the equivalence with connections can be found in \cite{SW1}. In the version presented here we give a global proof exploiting properties of flows of vector fields on compact manifolds. The idea is to think of a connection on a vector bundle as a compatible way of lifting vector fields $X$ on the base manifold to vector fields on the total space of the bundle, or equivalently to $X$-derivations acting on sections of the bundle. Even better, a connection can be described as a compatible way of lifting $\bR$-actions (determined by the flows of vector fields) on the base  manifold to $\bR$-actions on the covering bundle. 

Aside from its classical flavor, the equivalence of connections and parallel transport allows a description \cite{DST} of one dimensional topological field theories over a manifold (see Atiyah \cite{A} for a definition of topological field theories in arbitrary dimensions). Two dimensional topological field theories over a manifold $M$ admit a similar description in terms of connections on vector bundles over $LM$ the free loop space of $M$, along with some Frobenius data between fibers encoding parallel transport along pairs of pants cf. \cite{D3}. A careful extension of the properties of flows of vector fields used here to the world of supermanifolds leads us to a characterization of supersymmetric one dimensional topological field theories over a manifold cf. \cite{D2}.

\section{Preliminaries on flows.}
In what follows, $M$ denotes a compact manifold.

\begin{prop} \label{sum} Let $X$ and $Y$ be vector fields on $M$, and let $\a, \b:\bR\times M\ra M$ denote the flows determined by $X$, respectively $Y$. Then the flow $\g$ of the vector field $X+Y$ is given by 
\[ \g_t(x)= \lim_{n\to\infty}  \underbrace{(\a_{\frac{t}{n}}\b_{\frac{t}{n}})\comp \ldots \comp (\a_{\frac{t}{n}}\b_{\frac{t}{n}})}_{n}(x). \]

\end{prop}

\noindent\textit{Remark.} This is a version of the Trotter formula for flows on manifolds.\\

\noindent The proof is not hard: if we simply consider the composition of the flows, we get infinitesimally the vector field $X+Y$, but not a group action. The above zig-zag composition still 
generates infinitesimally $X+Y$ and defines a group action as well. Details can be found in \cite{D2}.

\begin{prop} \label{fv} Let $\a:\bR\times M\ra M$ be the flow of a vector field $X$ on the compact manifold $M$. If $f$ is a (positive) function on $M$ then the flow of $fX$ is given by $$\b:\bR\times M \ra M: (t,x)\mapsto \a(s(t,x),x),$$ where $s:\bR\times M \ra \bR$ is the solution to

$$\left\{ \begin{array}{l} 
\frac{\partial s}{\partial t}(t,x)= f(\a(s(t,x),x))
\\\\
s(0,x)=0, \text{ for all } x.

\end{array} \right.$$

\end{prop}

\noindent The proof is a routine check.

\begin{cor}  Let $X$ and $Y$ be vector fields on $M$. Then $X$ and $Y$ have the same (directed) trajectories if and only if $Y=fX$, for some positive function $f$ on $M$. 

\end{cor}

\begin{cor} If $Y=fX$, where $f$ is a positive function  on $M$, and $c$ is an integral curve of $X$ then $c\comp \f$ is an integral curve of $Y$, for some (orientation-preserving) diffeomorphism $\f$ of $\bR$.
\end{cor}

\section{Main Result}
After the above preliminaries we can state the main theorem. Consider an $n$-dimensional vector bundle $E$ over a compact manifold $M$. Denote by $\sP(M)$ the pathspace of $M$, i.e. the space of all (piecewise) smooth paths $\g:I\ra M $, where $I$ denotes an arbitrary interval. Then

\begin{thm}  \label{0} There is a natural 1-1 correspondence
\[  \left\{
\begin{array}{l}
\text{   Connections}\\
\text{ on $E$ over $M$}
\end{array} \right\}
\longleftrightarrow \left\{
\begin{array}{l}
\text{ Parallel transport maps }\\
\text{  associated to  $E$ over $M$}
\end{array} 
\right\}. \]
\end{thm}

Recall (compare with \cite{D4}) that a parallel transport map $\sP$ on $E$ over $M$ is a (smooth) section of the pullback bundle
\[ \xymatrix @C=4pc{ p^*Hom(E,E) \ar[r] \ar[d] & Hom(E,E) \ar[d] \\
\sP(M) \ar[r]^{p:=(i,e)} & M\times M,
} \]
where $Hom(E,E)$ denotes the bundle whose fiber at $(x,y)\in M\times M$ is given by $Hom(E_x, E_y)$, and $i(\g)$, $e(\g)$ denote the starting respectively the ending point of the path $\g$. In other words, the map $\sP$
associates to each path $\g$ in $M$ a linear map $\sP(\g):E_{i(\g)}\ra E_{e(\g)}$. This correspondence is required to satisfy the following:
\begin{enumerate}
\item $\sP(\g_x)=1_{E_x}$, where $\g_x$ is a constant map at $x\in M$.
\item (Invariance under reparametrization) $\sP(\g\comp\a)=\sP(\g)$, where $\a$ is an orientation-preserving diffeomorphism of intervals.
\item (Compatibility under juxtaposition) $\sP(\g_2\star\g_1)= \sP(\g_2)\comp \sP(\g_1)$, where $\g_2\star\g_1$ is the juxtaposition of $\g_1$ and $\g_2$.
\end{enumerate}
\begin{proof}
It is well known how a connection $\na$ gives rise to parallel transport: given a path $c$ in $M$, pull-back the connection along $c$ and solve a differential equation $(c^*\na) s=0$. The solutions $s$ to this differential equation are the parallel sections along $c$ and will define an isomorphism between the fibers at the endpoints of the path; see for example \cite{KN1} for details. 

We now need to show how a parallel transport map $\sP$ associated to the bundle $E$ over $M$ gives rise to a connection on $E$ over $M$. For this, let $X$ be a vector field on $M$ and $\a: \bR\ra \text{Aut}(M)$ the flow of $X$, where $\diff(M)$ stands for the group of smooth diffeomorphisms of $M$. Denote by $P$ the frame bundle of $E$, i.e. $P= GL(E)$, and let $\pi:P\ra M$ be the projection map.  Let  $\diff^G(P)$ denote the group of $G$-equivariant diffeomorphisms of $P$ or, equivalently, the group of diffeomorphisms of $P$ that descend to diffeomorphisms of $M$. Parallel transport allows us to lift the flow $\a$ of $X$ to a group homomorphism $\ti{\a}: \bR \ra \diff^G(P)$ which is the flow of a vector field $\ti{X}$ on $P$ that is $G$-invariant. Moreover, the vector field $\ti{X}$ is an $X$-derivation, i.e.
\[ \ti{X}:\cinf(P)\ra \cinf(P),\ \ \ \ \ \ \ \ti{X}(\pi^*(f)g)= \pi^*(X(f))g+f\ti{X}(g), \]
for $f\in\cinf(M)$ and $g\in\cinf(P)$. This holds since $\ti{X}$ is $G$-invariant and therefore $\ti{X}(\pi^*(f))=\pi^*(X(f))$.

Further, $\ti{X}$ extends to an $X$-derivation on $\cinf(P;V)$, the space of $V$-valued functions on $P$, where $V$ is another notation for $\bR^n$. Again, since $\ti{X}$ is $G$-invariant, it preserves $\cinf(P;V)^G$ the space of $G$-equivariant $V$-valued functions on $P$, which canonically identifies with $\G(E)$, the space of sections on $E$. We have defined a correspondence
\[ Vect(M)\ni X \longmapsto \{\ti{X}:\G(E)\ra\G(E)\}. \]
This defines  a connection if the following two conditions hold:
\begin{enumerate}
\item \hspace{4cm} $\widetilde{X+Y}= \ti{X}+\ti{Y}$
\item \hspace{4cm} $\widetilde{fX}=f\ti{X}$,
\end{enumerate}
for $X, Y$  vector fields on $M$, and $f$ a smooth function on $M$. Let us first show property (1). Consider $X, Y$ vector fields on $M$, with flows $\a$, respectively $\b$, and let $\g$ denote the flow of $X+Y$. Then 
\begin{eqnarray*}
\widetilde{X+Y} &=& \frac{d}{dt} \Big | _{t=0}\  \ti{\g}_t   \\ &=& \frac{d}{dt} \Big | _{t=0} \ \widetilde{ \lim_{n \to\infty}(\a_{t/n}\b_{t/n})^{(n)} }  \\
&=& \frac{d}{dt} \Big | _{t=0} \ \lim_{n\to\infty} (\ti{\a}_{t/n}\ti{\b}_{t/n})^{(n)} \\  &=& \ti{X}+\ti{Y}.
\end{eqnarray*} 
The third equality follows from the compatibility of lifting paths via parallel transport with concatenation of paths. For the second property (2), let $f$ be a (positive) function on $M$ and $X$ a vector field on $M$, with flow $\a$. By the Proposition \ref{fv} above, the flow $\a^f$ of the vector field $fX$ is given by the composition
\[ \xymatrix{ \bR\times M \ar[r]^-{1\times \D} & \bR\times M\times M \ar[r]^-{s\times 1} & \bR\times M \ar[r]^-{\a} & M,} \]
for some $s:\bR\times M\ra \bR$. The flow $\widetilde{\a^f}$ of $\widetilde{fX}$ is then given by the composition
\[ \xymatrix{ \bR\times P \ar[r]^-{1\times \D} & \bR\times P\times P \ar[r]^-{\ti{s}\times 1} & \bR\times P \ar[r]^-{\ti{\a}} & P,} \]
where $\ti{s}:\bR\times P \ra \bR$ is defined by the composition
\[ \xymatrix{ \bR\times P \ar[r]^-{1\times \pi} & \bR\times M \ar[r]^-{s} & \bR.} \]
But the above composition $\ti{a}(\ti{s}\times 1)(1\times\D)$ is nothing else than the flow of $f\ti{X}$. The vector fields $\widetilde{fX}$ and $f\ti{X}$ have the same flow, so they must be the same. The connection in the direction of the vector field $X$ is now determined by defining
\[ \na_X:= \ti{X}: \G(E)\ra\G(E).\]
It is clear that the family of derivations $\{\na_X\}$ parametrized by the space $Vect(M)$ of  vector fields on the manifold $M$ defines a connection $\na$ on the bundle $E$ over $M$. 

We are left to check that the two constructions are inverse to each other. The map ``$\longrightarrow$" which associates to a connection its parallel transport map is injective, as it is well known that the parallel transport of a connection {\it recovers} the connection.

The proof is complete if we can show that 
\[ \longrightarrow\comp\longleftarrow \ = id. \]
Indeed, start with a parallel transport map $\sP$ associated to the bundle $E$ over $M$. Denote by $\rho$ the standard representation of $G$ on $V=\bR^n$. Then $E$ can be canonically identified with the associated bundle $P\times_\rho V$, whose elements are classes $[p,v]\in P\times V/\sim$ (we say that $(p,v)$ and $(p',v')$ are in the same class if $p'=pg$ and $v'=g^{-1}v$, for some $g\in G$). Recall the isomorphism
\[ \cinf(P;V)^G\longrightarrow \G(E): \ \ \ \ \ \  f\longmapsto \{s(x)= [p,f(p)],\ \ \text{where } \pi(p)=x\}. \]
The parallel transport then defines as above a connection $\na$ given by
\[ \na_X[p,f(p)]= [p,\ti{X}(f)(p)],  \] 
in the direction of a vector field $X$ on $M$ ($\ti{X}$ denotes the lift of $X$ to the bundle $P$ via parallel transport as well as  the extension as a derivation to $V$-valued functions on $P$). Now, if $c:\bR\ra M$ is a curve in $M$ and $f\in\cinf(c^*P;V)^G$, we also have
\[ (c^*\na)_{\partial_t}[p,f(p)]= [p,v(f)(p)], \] 
where $v$ denotes the lift of $\partial_t$, the standard vector field on $\bR$, to $c^*P$. Therefore, a section $s=[p,f(p)]$ along $c$ is parallel with respect to the connection $\na$ if and only if $f$ is constant in the direction of the vector field $v$, which happens if and only if $s$ is parallel along $c$ with respect to the parallel transport $\sP$. This finishes the proof of the theorem.
\end{proof}

%Denote by  by $\text{Aut}(P,M)$ the group of automorphisms of $P$ that descend to $M$, i.e. that respect the fibers of $P$ over $M$. 
Let $q:\diff^G(P)\ra \diff(M)$ denote the obvious descending map. In the proof of Theorem \ref{0} we only used parallel transport along flows  of vector fields. This allows us to redefine a connection (on the principal $G$-bundle $P$) as a lift
\[  \left\{
\begin{array}{l}
\text{  Homomorphisms}\\
\ \ \ \bR \ra \diff(M)
\end{array} \right\}
\longrightarrow \left\{
\begin{array}{l}
\text{ Homomorphims }\\
\ \bR \ra \diff^G(P)
\end{array} 
\right\} \]

\[ \a \ \ \ \ \ \ \ \ \ \  \longmapsto \ \ \ \ \ \ \ \ \ \ \ \ti{\a} \]

\noindent that preserves the zig-zag composition of Proposition \ref{sum} and the action of the space of (positive) functions on $M$ given by Proposition \ref{fv}. Here, ``lift" means that a homomorphism $\a$ must map to a homomorphism $\ti{\a}$ that descends to $\a$, i.e. such that $q\ti{\a}=\a$. The lift should also preserve the trivial homomorphisms.

It is an interesting problem to see how  the {\it flatness} condition for a connection can be expressed in terms of lifting the dynamics of a manifold as above. It should be a lift that preserves the flow formula for the Lie bracket, but we are not aware of such a formula.\\

\noindent\emph{Acknowledgements.}  This note is based on an idea of Stephan Stolz that connections allow liftings of group actions. We warmly acknowledge the financial support provided by a BitDefender scholarship while visiting the Institute of Mathematics of the Romanian Academy ``Simion Stoilow".

\bibliographystyle{plain}
\bibliography{bibliografie}

\bigskip
\raggedright Institute of Mathematics of the Romanian Academy ``Simion Stoilow"\\  21 Calea Grivitei \\
Bucharest, Romania.\\ Email: {\tt florinndo@gmail.com}

\end{document}